# TOTAL POSITIVITY IN COORDINATE ALGEBRAS

G. Lusztig

## Introduction

**0.1.** Let $I$ be a finite set with a given symmetric bilinear form $\nu, \nu' \mapsto \nu : \nu'$ on $\mathbf{Z}[I]$ with values in $\mathbf{Z}$ such that $i : i = 2$ for any $i \in I$ and $i : j \in \{0, -1\}$ for any $i \neq j$ in $I$. (This is a Cartan datum of simply laced type.) Let $\mathbf{u}$ be the associative $\mathbf{C}$-algebra with generators $e_i (i \in I)$ and Serre relations $e_i^2 e_j - 2e_i e_j e_i + e_j e_i^2 = 0$ for any $i \neq j$ in $I$ such that $i : j = -1$, $e_i e_j - e_j e_i = 0$ for any $i \neq j$ in $I$ such that $i : j = 0$. We have $\mathbf{u} = \oplus_{\nu \in \mathbf{N}[I]} \mathbf{u}_\nu$ where the subspaces $\mathbf{u}_\nu$ are characterized by the properties that $\mathbf{u}_\nu \mathbf{u}_{\nu'} \subset \mathbf{u}_{\nu + \nu'}$, $1 \in \mathbf{u}_0$, $e_i \in \mathbf{u}_i$. We have $\dim \mathbf{u}_\nu < \infty$ for any $\nu$. Now $\mathbf{u}$ is a Hopf algebra with comultiplication $\Delta : \mathbf{u} \to \mathbf{u} \otimes \mathbf{u}$ characterized by $e_i \mapsto e_i \otimes 1 + 1 \otimes e_i$ for $i \in I$.

Let $\hat{\mathbf{u}} = \prod_{\nu \in \mathbf{N}[I]} \mathbf{u}_\nu$. This can be viewed as a completion of $\mathbf{u}$; the multiplication and comultiplication of $\mathbf{u}$ extend by continuity to a multiplication and comultiplication of $\hat{\mathbf{u}}$ making it into a Hopf algebra.

**0.2.** Let $\mathcal{O}$ be the set of $\mathbf{C}$-linear functions $\hat{\mathbf{u}} \to \mathbf{C}$ which are zero on all but finitely many $\mathbf{u}_\nu$. This is a commutative $\mathbf{C}$-algebra (a coordinate algebra) in which the product of $f', f''$ in $\mathcal{O}$ is the function $f'f'' : x \mapsto \sum_h f'(x'_h) f''(x''_h)$ where $\Delta(x) = \sum_h x'_h \otimes x''_h$.

**0.3.** Let $U(\mathbf{R}_{>0})$ be the totally positive semigroup associated in [L20,2.2] to our Cartan matrix and to $\mathbf{R}_{>0}$. This is the semigroup with generators $i^a, i \in I, a \in \mathbf{R}_{>0}$ and relations

$i^a i^b = i^{a+b}$ for $i \in I$ and $a, b$ in $\mathbf{R}_{>0}$

$i^a j^b i^c = j^{bc/(a+c} i^{a+c} j^{ab/(a+c)}$ for $i, j$ in $I$, $i : j = -1$ and $a, b, c$ in $\mathbf{R}_{>0}$

$i^a j^b = j^b i^a$ for $i, j$ in $I$, $i : j = 0$ and $a, b$ in $\mathbf{R}_{>0}$.

The following result will be verified in §1.

**Lemma 0.4.** *(a) There is a well defined map $\iota : U(\mathbf{R}_{>0}) \to \hat{\mathbf{u}}$ such that*

$\iota(i^a) = \sum_{p \in \mathbf{N}} a^p e_i^{(p)}$ *for any $i \in I, a \in \mathbf{R}_{>0}$ (we write $e_i^{(p)}$ instead of $e_i^p/p!$),*

$\iota(xy) = \iota(x)\iota(y)$ *for any $x, y$ in $U(\mathbf{R}_{>0})$.*

*(b) It satisfies*

$\Delta(\iota(x)) = \iota(x) \otimes \iota(x)$

*for any $x \in U(\mathbf{R}_{>0})$.*

**0.5.** Let $i_* = (i_1.i_2, \dots, i_m)$ be a sequence in $I$. For any $f \in \mathcal{O}$ and any $(a_1, a_2, \dots, a_k) \in \mathbf{R}^k_{>0}$ we have

$$\begin{aligned} &f(\iota(i_1^{a_1})\iota(i_2^{a_2})\dots\iota(i_k^{a_k})) \\ &= f(\sum_{(p_1,p_2,\dots,p_k)\in\mathbf{N}^k} a_1^{p_1}a_2^{p_2}\dots a_k^{p_k} e_{i_1}^{(p_1)}e_{i_2}^{(p_2)}\dots e_{i_k}^{(p_k)}) \\ &= \sum_{(p_1,p_2,\dots,p_k)\in\mathbf{N}^k} a_1^{p_1}a_2^{p_2}\dots a_k^{p_k} f(e_{i_1}^{(p_1)}e_{i_2}^{(p_2)}\dots e_{i_k}^{(p_k)}) \end{aligned}$$

where $f(e_{i_1}^{(p_1)}e_{i_2}^{(p_2)}\dots e_{i_k}^{(p_k)}) = 0$ for all but finitely many $p_1, p_2, \dots, p_k$. Thus,

$$(a_1, a_2, \dots, a_k) \mapsto f(\iota(i_1^{a_1})\iota(i_2^{a_2})\dots\iota(i_k^{a_k}))$$

can be regarded as a polynomial in $a_1, a_2, \dots, a_k$. We say that $f$ is in $\mathcal{O}_{\geq 0}$ if this polynomial has coefficients in $\mathbf{R}_{\geq 0}$ for any $i_*$ as above that is if

$$f(e_{i_1}^{(p_1)}e_{i_2}^{(p_2)}\dots e_{i_k}^{(p_k)}) \in \mathbf{R}_{\geq 0}$$

for any $i_*$ and any $(p_1, p_2, \dots, p_k)$ as above. Note that $\mathcal{O}_{\geq 0}$ is closed under addition and under scalar multiplication by elements of $\mathbf{R}_{\geq 0}$. By 0.4(b) it is also closed under multiplication. Thus $\mathcal{O}_{\geq 0}$ is a semiring. We say that $\mathcal{O}_{\geq 0}$ is the totally positive part of $\mathcal{O}$.

**0.6.** Let $\mathbf{B}_1$ be the canonical basis of $\mathbf{u}$ (a specialization at $v = 1$ of the canonical basis $\mathbf{B}$ of the quantum $v$-analogue $\mathbf{u}_v$ of $\mathbf{u}$). ($\mathbf{B}$ was defined in [L90] assuming that the Cartan datum is of finite type and in [L91], [K91] without this assumption. The definition in [K91] uses several ideas of [L90]; also, [K91] does not contain the positivity property [L91,11.5(a)] which is crucial for the present paper.)

For any $\mathbf{b} \in \mathbf{B}_1$ we define $\mathbf{b}^*$ to be the unique linear form $\hat{\mathbf{u}} \to \mathbf{C}$ which belongs to $\mathcal{O}$, takes $\mathbf{b}$ to 1 and takes any $\mathbf{b}' \in \mathbf{B}_1 - \{\mathbf{b}\}$ to 0. Note that $\mathbf{B}_1^* := \{\mathbf{b}^*; \mathbf{b} \in \mathbf{B}_1\}$ is a $\mathbf{C}$-basis of $\mathcal{O}$. Let $\mathcal{O}'_{\geq 0}$ be the subset of $\mathcal{O}$ formed by $\mathbf{R}_{\geq 0}$-linear combinations of elements of $\mathbf{B}_1^*$. The following result is proved in §2.

**Proposition 0.7.** *We have* $\mathcal{O}'_{\geq 0} \subset \mathcal{O}_{\geq 0}$.

**0.8.** We conjecture that

(a) the inclusion in 0.7 is an equality.

In type $A_2$ this follows from [L23,A4].

Assume that we are given a $\mathbf{C}$-vector space $V$ with a given basis $\beta$. We define $V_{\geq 0}$ to be the set of all linear combinations of elements of $\beta$ in which all coefficients are in $\mathbf{R}_{\geq 0}$. It is easy to see that $\beta$ is completelely determined by the set $V_{\geq 0}$ up to multiplication by elements in $\mathbf{R}_{>0}$.

It follows that if (a) holds then $\mathbf{B}_1$ can be reconstructed from $U(\mathbf{R}_{>0})$ (up to scalars in $\mathbf{R}_{>0}$).

One can hope that, similarly, the canonical basis [L93] of the (modified) full enveloping algebra can be reconstructed (up to scalars in $\mathbf{R}_{>0}$) from a positive version of a coordinate algebra. (See [L23,A5] for an example of this in type $A_1$.)

**0.9.** When the Cartan datum is no longer assumed to be simply laced,

$$\mathbf{u}, \hat{\mathbf{u}}, \mathcal{O}, U(\mathbf{R}_{>0}), \mathcal{O}_{\geq 0}$$

are still defined. But the the canonical basis is only defined (by the geometric definition in [L93]) up to sign. We believe that if we require that an element of the form $\mathbf{b}^*$ (of the dual canonical basis) is contained in $\mathcal{O}_{\geq 0}$ then the sign indeterminacy disappears. This should provide a simpler approach to the "normalization of signs" in [L93,19.2].

**0.10.** I thank Jeff York Ye for his comments on an earlier version.

## 1. Proof of Lemma 0.4

**1.1.** For $i \in I, a \in \mathbf{R}_{>0}$ we have

$$\begin{aligned}
\iota(i^a)\iota(i^b) &= \sum_p a^p e_i^{(p)} \sum_q b^q e_i^{(q)} \\
&= \sum_r \sum_{p,q;p+q=r} a^p b^q e_i^{(p)} e_i^{(q)} = \sum_r \sum_{p,q;p+q=r} a^p b^q \binom{r}{p} e_i^{(r)} \\
&= \sum_r (a+b)^r e_i^{(r)} = \iota(i^{a+b})
\end{aligned}$$

Now assume that $i, j$ in $I$, $i : j = -1$ and $a, b, c, A, B, C$ are in $\mathbf{R}_{>0}$. We set $\gamma = e_i e_j - e_j e_i \in \mathbf{u}$ and $\gamma^{(m)} = \gamma^m/m!$. Using [L93,42.1.2(f)] we have

$$\begin{aligned}
\iota(i^a)\iota(j^b)\iota(i^c) &= \sum_{p,q,r} a^p b^q c^r e_i^{(p)} e_j^{(q)} e_i^{(r)} \\
&= \sum_{p,q,r,n;n \leq p} \binom{p-n+r}{p-n} a^p b^q c^r e_i^{(q-n)} \gamma^{(n)} e_j^{(p-n+r)}.
\end{aligned}$$

Define $w, u$ in $\mathbf{N}$ by $p = n + w, q = n + u$. We obtain

$$\iota(i^a)\iota(j^b)\iota(i^c) = \sum_{u,w,r,n} \binom{w+r}{w} a^{n+w} b^{n+u} c^r e_i^{(u)} \gamma^{(n)} e_j^{(w+r)}.$$

We set $z = w + r$; we obtain

$$\begin{aligned}
\iota(i^a)\iota(j^b)\iota(i^c) &= \sum_{u,n,z} ( \sum_{w,r;w+r=z} \binom{z}{w} a^w c^r ) a^n b^{n+u} e_i^{(u)} \gamma^{(n)} e_j^{(z)} \\
&= \sum_{u,n,z} a^n (a+c)^z b^{n+u} e_i^{(u)} \gamma^n e_j^{(z)}.
\end{aligned}$$

Using [L93,42.1.2(g)] we have

$$\iota(j^A)\iota(i^B)\iota(j^C)=\sum_{p,q,r}A^pB^qC^re_j^{(p)}e_i^{(q)}e_j^{(r)}$$
$$==\sum_{p,q,r,n\le r}\binom{r-n+p}{r-n}A^pB^qC^re_i^{(r-n+p)}\gamma^{(n)}e_j^{(q-n)}.$$

Define $w,z$ in $\mathbf{N}$ by $r=n+w,q=n+z$. We obtain

$$\iota(j^A)\iota(i^B)\iota(j^C)=\sum_{z,w,p,n}\binom{w+p}{w}A^pB^{n+z}C^{n+w}e_i^{(w+p)}\gamma^ne_j^{(z)}.$$

We set $u=w+p$; we obtain

$$\iota(j^A)\iota(i^B)\iota(j^C)=\sum_{z,u,n}(\sum_{w,p;w+p=u}\binom{u}{w}A^pC^w)B^{n+z}C^ne_i^{(u)}\gamma^{(n)}e_j^{(z)}$$
$$=\sum_{z,u,n}(A+C)^uB^{n+z}C^ne_i^{(u)}\gamma^{(n)}e_j^{(z)}.$$

Assume now that $A=bc(a+c)^{-1},B=a+c,C=ab(a+c)^{-1}$. Then
$(A+C)^uB^{n+z}C^n=a^n(a+c)^zb^{n+u}$ for any $u,n,z$
hence

$$\iota(i^a)\iota(j^b)\iota(i^c)=\iota(j^A)\iota(i^B)\iota(j^C).$$

Next we assume that $i,j$ in $I$, $i:j=$ and $a,b$ are in $\mathbf{R}_{>0}$. We have

$$\iota(i^a)\iota(j^b)=\sum_{p,q}a^pb^qe_i^{(p)}e_j^{(q)}=\sum_{p,q}b^qa^pe_j^{(q)}e_i^{(p)}=\iota(j^b)\iota(i^a).$$

This completes the proof of Lemma 0.4(a)

We prove 0.4(b). Using (a) and the fact that $\Delta$ is an algebra homomorphism we see that it is enough to show that for $i\in I,a\in\mathbf{R}_{>0}$ we have $\Delta(\iota(i^a))=\iota(i^a)\otimes\iota(i^a)$. We have

$$\Delta(\iota(i^a))=\sum_pa^p\Delta(e_i^{(p)}=\sum_pa^p\sum_{q,r;q+r=p}e_i^{(q)}\otimes e_i^{(r)}$$
$$=\sum_{q,r;q+r=p}a^{q+r}e_i^{(q)}\otimes e_i^{(r)}=\iota(i^a)\otimes\iota(i^a).$$

This completes the proof of Lemma 0.4.

## 2. Proof of Proposition 0.7

**2.1.** It is enough to prove that for any $\mathbf{b} \in \mathbf{B}_1$ we have $\mathbf{b}^* \in \mathcal{O}_{\geq 0}$ that is

(a) $$\mathbf{b}^*(e_{i_1}^{(p_1)} e_{i_2}^{(p_2)} \dots e_{i_k}^{(p_k)}) \in \mathbf{R}_{\geq 0}$$

for any $i_1, i_2, \dots, i_k$ in $I$ and any $(p_1, p_2, \dots, p_k) \in \mathbf{N}^k$. But (a) is the coefficient of $\mathbf{b}$ when

$$e_{i_1}^{(p_1)} e_{i_2}^{(p_2)} \dots e_{i_k}^{(p_k)}$$

is written as a linear combination of elements of $\mathbf{B}_1$. This coefficient is in $\mathbf{N}$ by [L91,11.5(a)]. This proves (a). The proposition is proved.

Department of Mathematics, M.I.T., Cambridge, MA 02139